\magnification=1200
\def\za{\vrule height 6pt width 4pt depth 1pt}
\input amssym.tex
\font\eightrm=cmr8
\font\svntnrm=cmr17
\parskip=0pc
\baselineskip=14truept
\centerline{\svntnrm ON SOME EXAMPLES OF}
\bigskip
\centerline{\svntnrm  GROUP ACTIONS AND GROUP EXTENSIONS}
\vskip .5in
\centerline{\bf Alejandro Adem\footnote{*}{\eightrm{Partially supported
by an NSF grant, the MPIM-Bonn, Universit\'e de Paris 7 and CRM-Barcelona.}} 
and Erg\"un Yal\c{c}\i{}n}
\bigskip

\centerline{Department of Mathematics}

\centerline{University of Wisconsin}

\centerline{Madison, WI 53706 USA}
\vskip .5in

\noindent{\bf \S 0 Introduction}
\medskip
It is well-known that finite groups which can act freely on a sphere must be
rather special. To be precise we know that their abelian subgroups are cyclic
(see [CE], Theorem XII.11.6) and conversely that any group satisfying this
condition acts freely on a finite complex with the homotopy type of a sphere
(a result due to Swan [Sw]). How does this extend to other groups? The first fact
to keep in mind is that {\sl any} finite group will act freely on some product of
equidimensional spheres (see [O]). Hence it makes sense 
to determine the minimal number 
$n$ such that $G$ acts freely on a finite complex $X\simeq (\Bbb S^m)^n$
for some $m$. In this paper we address this problem for finite
$2$-groups; in this context the following interesting question has been
raised (see [ACPW], Benson's problem list):
\medskip
\noindent{\bf Question}: {\sl If $G$ is a finite $2$-group, is
$n$ equal to the rank of its largest elementary abelian subgroup?}
\medskip
\noindent Although there exist examples of groups of composite order for which
the answer is negative (the alternating group ${\cal A}_4$ is such
an example, see [O]), the methods used do not apply to the case
of $2$-groups, and this problem remains open even in the rank two situation.

In this note we present an analysis of this problem for certain classes of
$2$-groups, and show that it relates 
to interesting questions about group
extensions.
We first consider $2$--groups $G$ such that every
element of order $2$ is central (the 2C condition). 
If such a group acts freely on a
product $X\simeq (\Bbb S^m)^n$ we show that the action is homologically
trivial mod 2 and that the mod $2$ cohomology of the orbit space can be expressed
as $H^*(X/G,\Bbb F_2)\cong 
H^*(G, \Bbb F_2)/(\mu_1,\dots ,\mu_n)$ where $\mu_1,\dots ,\mu_n
\in H^{m+1}(G, \Bbb F_2)$ form a regular sequence of maximal length. After noting that
such an action can always be constructed using representation theory, we apply this
to provide a direct geometric proof of a recent cohomological characterization
of groups with the 2C property (see [AK]). 

Next we consider actions of groups expressed as an extension
$$1\to V\to G\to W\to 1\eqno (E)$$
where $V$ and $W$ are elementary abelian $2$--groups and $V$
is of maximal rank. Using cohomological arguments we find a bound on
the rank of $W$; in particular if the induced action of $W$ on $V$
permutes a basis, then we obtain
$$dim~W\le dim~V + dim V^W.$$
These methods are then used to show that if
$W=(\Bbb Z/2)^r$ acts freely on a finite complex $X$ with
cohomology generated by one--dimensional classes and such that the
action permutes a basis for $H_1(X,\Bbb F_2)$, then
$$r\le dim~H_1(X,\Bbb F_2) + dim~H_1(X,\Bbb F_2)^W
\le 2~dim~H_1(X,\Bbb F_2).$$
 From this we derive another 
\medskip
\noindent{\bf Question}: {If $W=(\Bbb Z/2)^r$ acts freely on a finite complex $X$ with cohomology
generated by one dimensional classes, does it follow that $r\le 2~dim~H_1(X,\Bbb F_2)$?}
\medskip
Using geometric restrictions,
we show that if a group $G$ defined as in (E) acts freely on 
$X\simeq (\Bbb S^m)^n$, where $n$ is the rank of $V$ and such
that $rk~W\ge 2n +1$, then the orbit space {\sl is not}
homotopy equivalent to $(\Bbb RP^m)^n$; indeed this indicates
that `exotic' intermediate orbit spaces must arise for free actions
of group extensions such as (E) where $|W|$ is large relative to $|V|$. 
However, for $(\Bbb RP^m)^n$ we
can show that if $(\Bbb Z/2)^r$ acts freely on it, then $r\le 2n$.

By using extensions of type (E) having $W$ of very
large rank (introduced by Ol'shanskii [Ol]), 
we construct some striking examples of very large 
elementary abelian $2$-groups
acting on products of projective spaces with relatively
small isotropy.
 For example, we show that there is an action of $(\Bbb Z/2)^{1249}$
on $X=(\Bbb RP^{2^{1298}-1})^{50}$ with isotropy subgroups of rank at
most $50$. Finally we use these groups to prove the following
\medskip
\noindent{\bf Proposition}

{\sl At least one of the two questions above has a negative answer.}~~~\za
\medskip

Unless stated otherwise,
all coefficients will be taken in $\Bbb F_2$, the field with two
elements, so they are suppressed throughout. In addition we will
assume that the dimensions of the spheres appearing in our arguments
are larger than $7$. We are grateful to 
Jon Alperin and George Glauberman for pointing out Ol'shanskii's
examples to us, to Dave Benson for his very helpful comments and to
the referee for his valuable advice.

\bigskip
\noindent{\bf \S 1 Free Group Actions and Orbit Spaces}
\medskip
Let $X$ denote a finite complex with a free action of a finite group $G$.
To compute the cohomology of the orbit space $X/G$ there is a spectral
sequence (Cartan-Leray) of the form
$E_2^{p,q}=H^p(G,H^q(X))$
converging to $H^{p+q}(X/G)$. In this note we will be interested in 
$2$--group actions on products of equidimensional spheres
(or on finite complexes homotopy equivalent to them). 

To begin we recall a result due to G. Carlsson [C]
for actions of $G=(\Bbb Z/2)^n$:
\medskip
\noindent{\bf Proposition 1.1}

{\sl If $G=(\Bbb Z/2)^n$ acts freely on $X\simeq (\Bbb S^m)^n$, trivially in homology, then
there exist classes $\mu_1,\dots ,\mu_n\in H^{m+1}(G)$ such that
$H^*(X/G)\cong H^*(G)/(\mu_1,\dots ,\mu_n)$.} ~~~\za
\medskip
The proof requires showing that we may take the $\mu_i$ to be transgressions
of the cohomology generators for the product of spheres and that these classes
form a {\sl regular sequence} $\mu_1,\dots ,\mu_n$ in $H^*(G)$. More generally
given a homologically trivial free action of any finite group on a product
of equidimensional spheres, the cohomology of the orbit space can be expressed
as a quotient of the type in 1.1 if the transgressions of the generators on the
fiber form a regular sequence.

Now by the results in [AB2], we know that any free action of $(\Bbb Z/2)^n$ on 
$(\Bbb S^m)^n$ is homologically trivial if $m > 7$ (our standing 
assumption). 
With this additional information
we can extend 1.1 to a more general class of groups.
\medskip
\noindent{\bf Theorem 1.2}

{\sl Let $G$ denote a finite $2$--group such that every element of order 2 in $G$
is central. If $G$ acts freely on $X\simeq (\Bbb S^m)^n$, where $n$ is the rank
of the largest elementary abelian $2$-subgroup in $G$,
then 
\smallskip
\noindent (1) $G$ acts homologically trivially on $X$ and 
\smallskip
\noindent (2) the transgressions of the cohomology generators for the spheres
form a regular sequence $\zeta_1,\dots ,\zeta_n\in H^{m+1}(G)$ and so
$$H^*(G)/(\zeta_1,\dots ,\zeta_n)\cong H^*(X/G).$$}

\noindent{\bf Proof:}
Let us first assume homological triviality.
We will use the fact that a group as above has Cohen--Macaulay
cohomology (a result due to Duflot [D]). It suffices to show
that $H^*(G)$ is finitely generated as a module over the subalgebra generated
by the $\zeta_i$. However note that there is a unique maximal elementary
abelian subgroup $E\subset G$ which has rank $n$ and is of course central.
By Quillen's $F$--isomorphism theorem, it suffices to show that 
$res^G_E(\zeta_1),\dots , res^G_E(\zeta_n)$ form a regular sequence in 
$H^*(E)$. This fact follows from 1.1.

To prove homological triviality, we consider the exact sequence

$$H^m(X/E)\longrightarrow H^m(X)\buildrel{d_{m+1}}\over
\longrightarrow H^{m+1}(E)$$
which arises from the spectral sequence associated to the action
of $E$. Note that by our hypotheses the $E$ action is homologically trivial.
As $G$ centralizes $E$, $G/E$ will act on this sequence, trivially
on $H^*(E)$. Now $d_{m+1}$ has rank $n$, hence it is injective and so
we deduce that 
the action must
be trivial on the entire cohomology, completing the proof.$~~~\za$

\medskip
More generally what we have is that if a $2$-group of rank $n$ acts
freely on a product of $n$ equidimensional spheres, then the centralizers
of the elementary abelian subgroups of maximal rank fulfill the conditions
above. In particular if $H^*(G)$ is Cohen--Macaulay, then the associated
transgressions $\mu_1,\dots ,\mu_n\in H^{m+1}(G)$ will be a regular
sequence (note that in addition the centralizers mentioned before will
detect the cohomology). Hence a lower bound on $m+1$ is given by the
minimal dimension which allows for a collection of equidimensional
regular elements which has the largest length.

Note that for an arbitrary finite group satisfying
the 2C condition
it would suffice to require that the $G$ action have odd order
isotropy subgroups to obtain a version of 1.2. Also we should mention
that Duflot's result is proved using a particular action of $G$ on a product
of spheres. The following construction provides one. Let $x_1,\dots ,x_n$ be
a basis for $E$, the unique maximal elementary abelian subgroup
in $G$, a group satisfying the 2C condition. Let $\chi_i$ denote a non-trivial 1-dimensional
real representation for $C_i=<x_i>$; denote $V_i=Ind^G_{C_i}(\chi_i)$.
Let $G$ act diagonally on $X=S(V_1)\times\dots\times S(V_n)\cong
(\Bbb S^{|G|/2-1})^n$; then one can easily verify that
the action restricted to the central subgroup $E\subset G$ is
free. As $E$ contains every involution in $G$, all the
isotropy subgroups must be of odd order.

More generally if $V\cong (\Bbb Z/2)^k$ 
is a central subgroup in $G$ the
above construction can be used to provide an action of $G$ on
$X=(\Bbb S^{|G|/2-1})^k$ such that if $\zeta_1,\dots ,\zeta_k$ are
the associated transgressions, then they form a regular sequence 
and
$$H^*(G)/(\zeta_1,\dots ,\zeta_k)\cong H^*_G(X)$$
where $H^*_G(X)$ denotes the equivariant cohomology of $X$. If $k$ is
equal to the rank of $G$ at $p=2$, we recover the $2$-central case
above.

The structural results above can be used to construct some very special
cohomology classes for finite $2$-groups. If $G$ is a $2$-group such that
every element of order 2 is central, then we may write it as
$G\cong G_0\times (\Bbb Z/2)^n$, where $G_0$ is a group (also satisfying
the $2$-central condition) such that every
maximal subgroup has the same rank as $G_0$ itself.
\medskip
\noindent{\bf Proposition 1.3}

{\sl Suppose that a group $G_0$ as above acts freely on 
$X=(\Bbb S^m)^n$, where $n$ is the rank of
$G_0$. Let $\pi_{G_0}:X/G_0\to BG_0$ denote the usual classifying map.
If $x\in H^{mn}(G_0)$ is such that $\pi_{G_0}^*(x)\ne 0$, then 
there exists a non-zero class $w\in H^{mn}(G_0)$ such that $res^{G_0}_H(w)=0$
for all proper subgroups $H\subset G_0$.}
\medskip
\noindent{\bf Proof:}
To prove this we take $H\subset G_0$ a maximal subgroup and consider
the following commutative diagram

$$\matrix{ H^N(G_0)&\buildrel{\pi_{G_0}^*}\over\longrightarrow &H^N(X/G_0)\cr
\uparrow&&\uparrow\cr
H^N(H)&\buildrel{\pi_H^*}\over\longrightarrow &H^N(X/H)\cr}$$
where the two vertical arrows are transfer maps
and $N=mn$. 
Note that as $H$ still has maximal rank, our previous arguments
show that $\pi_H^*$ is onto. Furthermore the right hand vertical
map is surjective, as its image is the kernel of multiplication
by the one dimensional class corresponding to $H$; as $N$ is the 
top dimension this is the zero map. If we take the
class $x\in H^N(G_0)$, then there exists a $y\in H^N(H)$ 
such that $tr\cdot \pi_H^* (y) = \pi_{G_0}^*(x)$; hence
if we take $z=tr(y)$, then $z$ is a class which satisfies the
same hypothesis as $x$, but is a transfer image. Continuing like
this, using index 2 subgroups, we can construct a class $w\in H^N(G_0)$
which is in the image of the transfer from $H^*(E)$ ($E$ the unique 
maximal elementary abelian subgroup) which still satisfies the same
condition that $x$ does. 

Now let $K$ denote any maximal subgroup, we will use the double coset
formula  to compute $res^{G_0}_K(w)$. If $w=tr^{G_0}_E(u)$, then the following
terms will appear in this formula: $tr^K_{K\cap E^g}\cdot 
res^{E^g}_{E^g\cap K}(c_gu)$, where $c_g$ is the map induced by
conjugation by $g$ and the summation is taken over $\overline{g}\in K\backslash G_0/E$.
In our case $E$ is central, hence the conjugation action is trivial,
$E=E^g$, $E^g\cap K=E$ and $K\backslash G_0/E$ has $[G_0:K]$ distinct elements.
Hence we obtain
$$res^{G_0}_K(w)=[G_0:K]u=0$$ 
and the proof is complete.
$~~~\za$
\medskip
Note that for the case of a general $G\cong G_0\times (\Bbb Z/2)^n$, we can take
$\gamma$ to be the product of all one dimensional classes in
$H^*(G_0\times (\Bbb Z/2)^n)$ which involve non--zero terms from the
elementary abelian factor, and if $w$ is chosen as above, then
$w\gamma$ will be a non--zero class restricting trivially on all proper 
subgroups in $G$. A more complicated but purely algebraic construction of 
such classes was
provided in [AK]. From this it can be shown that the finite groups which
have Cohen--Macaulay mod 2 cohomology with undetectable elements are precisely
the $2$-groups such that every element of order two is central.
\bigskip
\noindent{\bf \S 2 Group Extensions}
\medskip
In this section we will consider actions of $2$-groups which can be described
as extensions of the form

$$1\to V\to G\to W\to 1\eqno (2.0)$$
where $V$ and $W$ are both elementary abelian of finite rank. In addition we
will assume that $V$ is of maximal rank in $G$. We begin with a 
purely group--theoretic fact.
\medskip
\noindent{\bf Proposition 2.1}

{\sl Let $G$ be an extension as above, where $W$ acts on a basis for $V$ via
permutations. Then 
$dim~W\le dim~V + dim~V^W.$}
\medskip
\noindent{\bf Proof:} Choose $H\subset W$ the subgroup of elements which act
trivially on $V$ and consider the five-term exact sequence associated
to the extension $1\to V\to \tilde H\to H\to 1$:

$$0\to H^1(H)\to H^1(\tilde H)\to H^1(V)\buildrel{d_2}\over\to H^2(H)\to
H^2(\tilde H).$$
The quotient group $W/H$ acts on this, trivially on $H^2(H)$, hence we
obtain that the coinvariants $H^1(V)_W$ have rank at least as large as
the dimension of $im~d_2$. We claim that $im~d_2$ has rank at least
$rk(H)/2$. Suppose otherwise, then the k--invariants of the extension
would form a collection of say $q$ quadratic polynomials in $H^*(H)$,
a polynomial ring on at least $2q+1$ variables. This implies that they
have a common zero (see [G], pages 11 and 18); therefore an extra $\Bbb Z/2$ summand would
split off
and hence $\tilde H$
would have rank larger than that of $V$, which is impossible. Hence
by dualizing we conclude that
$rk~H\le 2~dim~V^W$. Now $W/H$ acts faithfully on $V$, a permutation module
and so one can easily show that 
$rk~W/H\le dim~V - dim V^W$
from which we can then infer that
$$rk~W\le dim~V + dim~V^W.
~~~\za$$
\medskip
\noindent{\bf Corollary 2.2}

{\sl If all elements of order 2 in $G$ are central
then $rk~G/V\le 2~rk~V$.}$~~~\za$
\medskip
More generally note that if $Q\subset GL(V)$ is elementary abelian, then
$rk~Q\le (rk~V)^2/4$, and so the proof above can be modified
to yield the {\sl quadratic} bound
$$dim~W\le 2~dim~V^W + (dim~V)^2/4.$$
\medskip 

We now switch to topology to prove a result analogous to 2.1 for
group actions.
\medskip
\noindent{\bf Theorem 2.3}

{\sl Let $X$ denote a finite complex such that $H^*(X)$ is generated by 
one dimensional classes. If $W=(\Bbb Z/2)^r$ acts freely on $X$
inducing a permutation module on $H_1(X)$, then

$$r\le dim~H_1(X) + dim~H_1(X)^W\le 2~dim~H_1(X).$$}
\medskip
\noindent{\bf Proof:}
As before we choose $H\subset W$, the subgroup of elements acting trivially
on $H^1(X)$, and hence trivially on all the cohomology. Consider the transgression
map $d_2:H^1(X)\to H^2(H)$ arising from the spectral sequence for the action
of $H$ (as discussed in \S 1). Now we claim as before that $im~d_2$ has 
dimension at least $rk(H)/2$. Suppose otherwise, then we can find a cyclic subgroup
$C\subset H$ such that $res^H_C\cdot d_2=0$, and hence the corresponding differential
$d_2$ will be identically zero in the spectral sequence for the action of
$C$. As $H^*(X)$ is generated by one dimensional classes, this means that in fact
the spectral sequence collapses, and so the $E_{\infty}$ term is infinite dimensional.
However the spectral sequence converges to the cohomology of the finite complex $X/C$,
whence we obtain a contradiction.

The rest of the proof now follows directly from the arguments used in proving 2.1.~~~\za
\medskip
 From the topological point of view, the following question seems
very natural:
\medskip
\noindent{\bf Question A}

{\sl If $W=(\Bbb Z/2)^r$ acts freely on a finite complex $X$
with cohomology generated by one dimensional classes, does it
follow that
$r\le 2~dim~H_1(X)$?}
\medskip
Now we consider actions on products of spheres by
group extensions of the type defined in this
section. As before let us assume that $V$ has maximal rank $n$,
and suppose that $G$ acts freely on $X\simeq (\Bbb S^m)^n$. As we have seen,
$V$ acts homologically trivially on $H^*(X)$. More generally, we
know by a Hopf invariant argument (see [AB2])
that $G$ acts on $H^m(X)$ by permuting the canonical basis. Hence
we have a fundamental exact sequence of $G/V=W$ modules:

$$0\to H^m(X)\buildrel{d_{m+1}}\over\longrightarrow
H^{m+1}(V)\longrightarrow H^{m+1}(X/V)\to 0$$
where the module on the left is a {\sl permutation module}.
Or equivalently, $W$ permutes the k--invariants $\mu_1,\dots ,\mu_n$,
where 
$$H^*(X/V)\cong \Bbb F_2[x_1,\dots ,x_n]/(\mu_1,\dots ,\mu_n).$$

This provides a geometric constraint on the possible group extension. An
immediate consequence is the following.
\medskip
\noindent{\bf Theorem 2.4}

{\sl If $X/V\simeq (\Bbb R P^m)^n$, then
$rk~W\le rk~V^W + rk~V.$}
\medskip
\noindent{\bf Proof:} Our hypotheses imply that 
$$H^*(X/V)\cong \Bbb F_2[y_1,\dots ,y_n]/(y_1^{m+1},\dots ,y_n^{m+1})$$
for one dimensional classes $y_1,\dots ,y_n$. If we consider these k-invariants
as elements in $H^*(V)$, then they must span a $W$-submodule $M$. Write 
$m+1=2^rk$, where $k$ is an odd number. If $k>1$, then by the binomial
coefficient theorem it is easy to see that the vectors $y_1,\dots ,y_n$
must be permuted by the action of $W$, as otherwise their $m+1$ powers 
would not span a $W$-submodule (note that this argument does not require any
homotopy theory). If $k=1$ and so $m+1$ is a power
of two, then the iterated squaring map induces a $W$--isomorphism 
between the permutation module $M$ and $H^1(X/V)$.
Applying 2.3 completes
the proof.$~~~\za$
\medskip
\noindent{\bf Corollary 2.5}

{\sl If $W=(\Bbb Z/2)^r$ acts freely on $Y\simeq (\Bbb R P^m)^n$, then
$r\le 2n$.}$~~~\za$

\medskip
If $m=4k +1$, the proof can be modified to show that in fact
$r\le n$. The point is that in this case the ideal generated by the
transgressions is invariant under the action of the
Steenrod algebra. A result due to Carlsson [C] implies that
if there were more than $n$ variables, then the quadratic polynomials
would have a common zero. This observation is due to Cusick [Cu]. 
Recall (see [AB1]) that the {\sl free $2$-rank of symmetry} of a finite complex $X$ is
the rank of the largest elementary abelian $2$-group which acts freely on it.
Our results extend those of Cusick [Cu] to actions which are
not homologically trivial, and can be summarized in the following
\medskip
\noindent{\bf Theorem 2.6}

{\sl 
If $X=(\Bbb R P^m)^n$ and $\phi (X)$ denotes its free $2$-rank of symmetry, then
we have

$$\phi (X) = \cases{0&if~$m\equiv 0,2$ mod 4\cr
                    n&if~$m\equiv 1$ mod 4\cr
                    2n&if~$m\equiv 3$ mod 4.\cr}$$}
\medskip
\noindent{\bf Proof:} We have already established the upper bounds for
$m$ odd. For $m$ even we have $\chi (X)=1$, whence no non-trivial
finite group can act freely on $X$. To complete the proof it suffices to
to observe that the standard free actions of $\Bbb Z/4$ on $\Bbb S^{4q+1}$
and of the quaternion group $Q_8$ on $\Bbb S^{4q+3}$ give
rise to actions on real projective space whose products achieve the
desired upper bounds.~~~\za

\medskip

\noindent{\bf Corollary 2.7}

{\sl An extension $G$ of type (2.0) acts freely and homologically trivially
on $X\simeq (\Bbb S^m)^n$ 
such that $X/V\simeq (\Bbb R P^m)^n$ if and only if
every element of order 2 in $G$ is central.}
\medskip
\noindent{\bf Proof:} We have already seen that any $2$-group satisfying the
2C condition acts freely and homologically trivially
on a product of equidimensional spheres such that
the generators of the
maximal (central) $2$-torus act via multiplication by -1 on the factors.
Hence $X/V\cong (\Bbb R P^m)^n$ in this case. Conversely the hypotheses
imply that the $W$ action on $H^1(V)$ is trivial, hence $V$ is central
and of maximal rank, which implies the 2C condition.$~~~\za$
\medskip

We now exhibit some interesting examples of group extensions.
\medskip
\noindent {\bf Examples 2.8}

The following examples appear in [Ol]. Let $\Phi = \{\phi_1,\dots ,\phi_t\}$
denote a collection of skew--symmetric bilinear forms on an $n$--dimensional
$\Bbb F_2$--vector space $M$. 
A group $G_{\Phi}$ can be defined as follows. As generators
we will take elements

$$a_1,\dots ,a_n,~~~~~~~b_1,\dots ,b_t$$
subject to the relations 

$$a_i^2=b_j^2=[a_i,b_j]=[b_s,b_u]=1$$
and
$$[a_i,a_j]=b_1^{\phi_1(a_i,a_j)}\dots b_t^{\phi_t(a_i,a_j)}$$
for $i,j=1,2,\dots ,n$ and $s,u=1,\dots ,t$.
$G_{\Phi}$ is a class 2 nilpotent group, $B=<b_1,\dots ,b_t>$ is central
and $G_{\Phi}/B\cong M=<a_1B,\dots ,a_nB>$.

In [Ol] it is shown that given an integer $k>0$ such that 
$n<t(k-1)/2$, then we may choose $\phi_1,\dots ,\phi_t$
such that the elementary abelian subgroups in the associated
group $G_{\Phi}$ are of
rank at most $t+k-1$.

 For example, if $k=51$, $t=50$ and $n=1249$, then $G_{\Phi}$ fits into
an extension of the form 

$$1\to(\Bbb Z/2)^T\to G_{\Phi}\to (\Bbb Z/2)^N\to 1$$
where $T\le 100$, $N\ge 1199$ and $T$ is the rank of $G_{\Phi}$.

Using the construction described in \S 1, we obtain an action of 
$G_{\Phi}$ on $X=(\Bbb S^m)^{50}$, $m=2^{1298}-1$, such that the generators of $B$ act
via multiplication by -1 on the spheres and of course the maximal
rank of an isotropy subgroup is $T-50\le 50$. Dividing out by the $B$ action
we obtain an action of $(\Bbb Z/2)^{1249}$ on $(\Bbb R P^m)^{50}$ with isotropy
of rank at most $50$.

To tie things together we now introduce a well--known question on group
actions (see [ACPW], D. Benson's problem list).
\medskip
\noindent{\bf Question B}

{\sl If $G$ is any finite $2$-group of rank $l$, does there exist
a finite complex $X\simeq (\Bbb S^m)^l$ such that $G$ acts freely
on $X$?}
\medskip
\noindent{\bf Proposition 2.9}

{\sl Questions A and B cannot both have affirmative answers.}
\medskip
\noindent{\bf Proof:} Suppose that question B were true and take 
$G_{\Phi}$ as before. Then we would have a free action of 
$(\Bbb Z/2)^N$ on $X/(\Bbb Z/2)^T$, where $X\simeq (\Bbb S^m)^T$. As $T$ is the
rank of the group, a maximal elementary abelian subgroup $V$ acts homologically trivially on 
$X$, and as we have seen
$$H^*(X/V)\cong \Bbb F_2 [x_1,\dots ,x_T]/(\mu_1,\dots ,\mu_T)$$
where the $x_i$ are one dimensional and the $\mu_j$ form a regular
sequence in the polynomial algebra. Note however that $N\ge 1199$,
whereas $T\le 100$! This negates question A; indeed $G_{\Phi}$
must provide a counterexample to one of the two questions. Hence a proof
of one of them will disprove the other one. $~~~\za$
\medskip
A concrete result which we do have is the following.
\vfill\eject
\noindent{\bf Theorem 2.10}

{\sl If $G$ as in (2.0) acts freely on $X\simeq (\Bbb S^m)^n$, where
$n=rk~V$ and $rk~W\ge 2n+1$, then $X/V$ is not homotopy
equivalent to $(\Bbb R P^m)^n$.} $~~~\za$
\medskip
It should be noted that the dimension of the spheres involved in our
free actions must be sufficiently large to produce enough homology.
More precisely, suppose we know that an extension $G$ as in (2.0) acts
freely on $X\simeq (\Bbb S^m)^t$, where we assume $dim~V=rk~G=t$. As we know,
the cohomology of the orbit space $X/V$ will be of the form
$\Bbb F_2[x_1,\dots ,x_t]/(\mu_1,\dots ,\mu_t)$; note that the total
dimension of the cohomology is $(m+1)^t$. Using cohomological
exponents it is possible to show that 
$$dim~G/V\le (m+1)t.$$
We sketch the proof: there is an index $2^t$ subgroup $H\subset G$ which 
acts trivially on $H^*(X,\Bbb Z)$.
Hence $H/H\cap V \subset G/V$ acts freely
and preserving orientation on $X/H\cap V$. Moreover it acts
trivially on the torsion free part of $H^*(X/H\cap V)$.
Decomposing the top class in terms of $m$ dimensional
classes and using a theorem of
Browder [B], we can obtain $log_2~|H/H\cap V|\le dim~X/H\cap V$, and so we
conclude $dim~G/V-t\le dim H/H\cap V\le mt$, whence the result follows.

A conjecture due to Carlsson (see [AB1]) would imply that in
all cases
$$2^{dim~G/V}\le (m+1)^t.$$
 For example if we take $G_{\Phi}$ as before, then 
the first inequality yields
$1099\le 100m$, whence $m\ge 11$. The conjectured
inequality would yield
$2^{1199}\le (m+1)^{100}$ from which we would deduce 
$m\ge 2^{11}-1$.

\bigskip\bigskip
\centerline{\bf REFERENCES}
\medskip
\item{[AB1] } Adem, A. and Browder, W., ``The Free Rank of Symmetry
of $(\Bbb S^n)^k$,'' Inventiones Math. 92 (1988), 431--440.
\smallskip
\item{[AB2] } Adem, A. and Benson, D., ``Elementary Abelian Groups 
Acting on Products of Spheres,'' Math. Zeitschrift, to appear.
\smallskip
\item{[ACPW] } Adem, A., Carlson, J., Priddy, S., Webb, P.,
(editors) {\bf Proceedings of Symposia in Pure Mathematics},
Seattle 1996 AMS Summer Institute Volume (1997).
\smallskip
\item{[AK] } Adem, A. and Karagueuzian, D., ``Essential Cohomology
of Finite Groups,'' Comm. Math. Helvetici 72 (1997), 101--109.
\smallskip
\item{[B] } Browder, W., ``Cohomology and Group Actions,''
Inventiones Math. 71 (1983), 599--607.
\smallskip
\item{[C] } Carlsson, G., ``On the Nonexistence of Free Actions of
Elementary Abelian Groups on Products of Spheres,'' Am. J. Math.
102 (1980), 1147--1157.
\smallskip
\item{[CE] } Cartan, H. and Eilenberg, S., {\bf Homological Algebra},
Princeton University Press (1956).
\smallskip
\item{[Cu] } Cusick, L.. ``Elementary Abelian $2$-Groups that Act Freely
on Products of Real Projective Spaces,'' Proc. Am. Math. Soc. Vol.87,
No.4 (1983), 728--730.
\smallskip
\item{[D] } Duflot, J.,``Depth and Equivariant Cohomology,''
Comm. Math. Helv. 56 (1981), 627--637.
\smallskip
\item{[G] } Greenberg, M., {\bf Lectures on Forms in Many Variables},
Benjamin, New York 1969.
\smallskip
\item{[O] } Oliver, B., ``Free Compact Group Actions on Products of
Spheres,'' in Springer-Verlag Lecture Notes in Mathematics, Vol. 763
(Arhus 1978 Algebraic Topology Volume), pp. 539--548.
\smallskip
\item{[Ol] } Ol'shanskii, A., ``The Number of Generators and Orders
of Abelian Subgroups of Finite $p$--groups,'' Math. Notes Vol. 23.
No. 3 (1978), 183--185.
\smallskip
\item{[Sw] } Swan, R. G., ``Periodic Resolutions for Finite Groups,''
Annals Math. 94 (1960), 267--291.
\end